\documentclass[final]{article}

\usepackage{graphicx}
\usepackage{amsmath}
\usepackage{amsfonts}
\usepackage{amssymb}
 \usepackage{amscd}

\newtheorem{example}{Example}[section]

\newtheorem{remark}{Remark}[section]
\newtheorem{con}{Conjecture}[section]
\newtheorem{thm}{Theorem}[section]
\newtheorem{lem}{Lemma}[section]
\newtheorem{cor}{Corollary}[section]
\newtheorem{prop}{Proposition}[section]

\newenvironment{proof}[1][Proof]{\textbf{#1.} }{\ \rule{0.5em}{0.5em}}

\begin{document}

\title{ New log-majorization results concerning eigenvalues and singular values and a complement of a norm inequality}

\author{{\small Mohammad M. Ghabries}$^{\rm a}$\thanks{Email:mahdi.ghabries@gmail.com},  \vspace{3pt}
{\small Hassane Abbas}$^{\rm b}$, 
\vspace{2pt} {\small Bassam Mourad}$^{\rm b}$,\vspace{2pt} {\small Abdallah Assi}$^{\rm a}$ \vspace{2pt} \\
 $^{\rm a}${\small{LAREMA, Facult\'e des Sciences-D\'epartement de math\'ematiques, Angers, France}};\\
  $^{\rm b}${\small{Department of Mathematics, Faculty of Science, Lebanese University, Beirut, Lebanon}}}

\maketitle

\begin{abstract}
The purpose of this paper is to establish  new log-majorization results concerning eigenvalues and singular values which generalize some previous work related to a conjecture and  an open question which were presented by R. Lemos and G. Soares in \cite{lemos}. In addition, we present a complement of a unitarily invariant norm inequality
which was conjectured by R. Bhatia, Y. Lim and T. Yamazaki in \cite{Bhatia2}, and recently proved by T.H. Dinh, R. Dumitru and J.A. Franco in \cite{Dinh} for the Schatten p-norm with $1\leq p\leq \infty$.
\end{abstract}

\paragraph*{keywords.}{\footnotesize  Positive semi-definite matrix; Hermitian matrix; Log-majorization; Eigenvalues; Singular values; Norm Inequalities}


\section{Introduction}

Let $M_{n}$ be the space of $n\times n$ complex matrices where its identity matrix is denoted by $I_{n}$. If $A$ is positive semi-definite, then we will write
$A\geq 0$ to indicate so. If,  in addition,  $A$ is invertible, then we will write $A>0$.
For positive semi-definite matrices $A, B\in M_{n}$,  we write $A \geq B$ if $A - B\geq 0$. For any $n\times n$ matrix $X$, we shall assume that its singular values, which are the eigenvalues of the positive semi-definite matrix $\vert X\vert = (X^*X)^{1/2}$,
are taken in the decreasing order: $s_1(X)\geq s_2(X)\geq \dots \geq s_n(X)\geq 0$. In addition, if the eigenvalues  $\lambda_1(X), \lambda_2(X), \dots, \lambda_n(X)$ of  $X$ are real, then we will always assume that they are also arranged  in decreasing order, that is, $$\lambda_1(X) \geq \lambda_2(X) \geq \dots \geq \lambda_n(X).$$

For a positive semi-definite matrix $X\in M_{n}$, we shall denote $$\lambda(X) = \left(\lambda_1(X),\lambda_2(X),\dots,\lambda_n(X)\right)^{t}$$ which is clearly a real column vector, and so is the vector $s(X):=(s_1(X), s_2(X), \dots,s_n(X))^{t}$.

If $\lambda(A)$, $\lambda(B)$ $\in \mathbb{R}^n,$ then by $\lambda(A) \prec_{w log} \lambda(B)$, we mean that $\lambda(A)$ is \emph{weakly log-majorized} by $\lambda(B)$, that is  $$\prod_{i=1}^{k}\lambda_i(A) \leq \prod_{i=1}^{k}\lambda_i(B) \hspace{1cm} \mbox{ for all }  k = 1,2,\dots,n.$$  In addition, we shall write $\lambda(A) \prec_{log} \lambda(B)$ and we will say that $\lambda(A)$ is \emph{log-majorized} by $\lambda(B)$    if the preceding inequality is true for
$ k = 1,2,\dots,n-1$ and  equality holds for $k = n$ \  i.e.  $$\prod_{i=1}^{n}\lambda_i(A) = \prod_{i=1}^{n}\lambda_i(B).$$


We denote by $\vert\vert\vert \cdot\vert\vert\vert$ any unitarily invariant norm on the space $M_n$ and by $\vert\vert .\vert\vert_p$ for $1\leq p\leq \infty$ the Schatten p-norm, which is defined by
$$ \|X\|_{p}=\mathrm {Tr} (|X|^{p})^{\frac {1}{p}}  =\left(\sum\limits_{i=1}^{n}s_i^p(X)\right)^{1/p} .$$
Of special importance are the cases  $p = 1$ named the trace norm, $p = 2$ the Hilbert-Schmidt norm if $p < \infty$ and $\vert\vert X\vert\vert_{\infty} = s_1(X)$ which is known as the operator norm.\\

The weighted geometric mean of two positive definite matrices $A$ and $B$ is defined by $$A\sharp_{t}B = A^{\frac{1}{2}} (A^{-\frac{1}{2}} B A^{-\frac{1}{2}})^{t} A^{\frac{1}{2}} \hspace{1cm} \ \ \ \ \ \ \ \ \ \text{for} \ \  0\leq t\leq 1.$$ As it is well known, this can be extended  to the positive semi-definite matrices as  follows $$A\sharp_{t}B = \lim_{\epsilon \rightarrow 0^+}(A + \epsilon I_n)\sharp_t (B + \epsilon I_n) \hspace{1cm} \text{for} \ \ 0\leq t\leq 1.$$

One important fact is that $A\sharp_t B = B\sharp_{1 - t}A$ for all $A, B \geq 0$. Clearly, if $t = \frac{1}{2}$, then we obtain the well known geometric mean $A\sharp B = A^{\frac{1}{2}} (A^{-\frac{1}{2}} B A^{-\frac{1}{2}})^{\frac{1}{2}} A^{\frac{1}{2}}.$ \\

Let $A$ and $B$ be two $n\times n$ positive semi-definite matrices, L. Zou \cite{Zou} proved the following inequality \begin{equation} \prod_{j=1}^{k} s_j\left(A^{\frac{1}{2}} (A\sharp B) B^{\frac{1}{2}}\right) \leq \prod_{j=1}^{k} s_j(AB), \ \ \ \ \ \ k = 1, \dots, n.  \label{1}  \end{equation}

Recently, R. Lemos and G. Soares \cite{lemos} gave another proof of \eqref{1}, and asked whether it is possible to find some  generalization of it. More
explicitly, their question gives rise to  the following conjecture.

\begin{con}Let $A$ and $B$ be two positive semi-definite matrices. Then \begin{equation} s\left(A^t (A\sharp_t B) B^{1 - t}\right) \prec_{log} s(AB), \hspace{1cm} \ \ \ \ \ \ \ \ \ \ \ \ \ \ \  0\leq t\leq 1.\label{con}\end{equation} \end{con}

In the same paper they proved that for all $A, B\geq 0$ and $r, s \in \mathbb{R}$, \begin{equation}\lambda\left( (A\sharp_{r,t}B)(A\sharp_{s, 1 - t}B) \right) \prec_{log} \lambda( A^{r + s - 1}B) \ \ \text{for} \ \  0\leq t\leq 1,\label{2}\end{equation} where $$A\sharp_{r,t}B = A^{\frac{r}{2}} (A^{-\frac{1}{2}} B A^{-\frac{1}{2}})^{t} A^{\frac{r}{2}}.$$

Earlier in \cite{Hiai}, F. Hiai and M. Lin proved inequality \eqref{2} for $r = s = 1$ and asked whether it is possible to replace the eigenvalues with singular values. More generally, R. Lemos and G. Soares \cite{lemos} conjectured the following inequality.

\begin{con} If $A, B \geq 0, r, s \in \mathbb{R}$ and $0\leq t\leq 1$, then \begin{equation}s\left( (A\sharp_{r,t}B)(A\sharp_{s, 1 - t}B)\right) \prec_{log} s\left( A^{r + s - 1}B\right).\label{ss}\end{equation}
\end{con}

\begin{remark} The case  $r = s = 1$ was proved by F. Hiai and M. Lin precisely for $\frac{1}{4}\leq t\leq \frac{3}{4}$ in \cite{Hiai}. Later, R. Lemos and G. Soares \cite{lemos}  showed that Conjecture 1.2 is valid for $r, s \geq 0$ and $\frac{r}{r + s} \leq 2t \leq \frac{2r + s}{r + s}$. For the latter case, we
shall prove further generalization of \eqref{ss} (see Theorem 2.1 below).
\end{remark}

Our main objective in this paper is to present new log-majorization results related to the above two conjectures as well as to prove an interesting norm inequality. More explicitly, the rest of the paper is organized as follows. In Section 2, we first disprove Conjecture 1.2 in its current setting by providing a counterexample in the case $t = 0$ and $r, s$ being two real numbers such that $s - 1\leq 0$ and $2r + s - 1\geq 0$. Then, we give a further generalization of  \eqref{2} as well as of \eqref{ss} in the cases where its has been shown to be valid. In Section 3, new results related to Conjecture 1.1 are established, and a new proof for a refinement of inequality \eqref{1} is presented. In the last section, we present a complement of a unitarily invariant norm inequality which was conjectured by R. Bhatia, Y. Lim and T. Yamazaki \cite{Bhatia2} and recently proved by T.H. Dinh, R. Dumitru and J.A. Franco \cite{Dinh} for the Schatten p-norm with $1\leq p\leq \infty$.

\section{Results concerning Conjecture 1.2}
This section deals first with  presenting  a counterexample to Conjecture 1.2. Then some generalization results related to this conjecture are given. Finally,  based on our results, we conclude with formulating an alternative conjecture.

The starting point for us is the following lemma whose proof can be found in \cite{GAM}.
\begin{lem} Let $A$ be a positive definite matrix and $B$ be a Hermitian matrix. Then for all $p, q \geq 0$ $$\lambda(A^p B A^{-q} B) \succ_{log} \lambda(A^{p - q} B^2).$$ \end{lem}

As a result, we have the following example which shows that Conjecture 1.2 is not valid in its current setting.

\begin{example} Let $r, s$ be two real numbers such that $s - 1\leq 0$ and $2r + s - 1\geq 0$. Substituting these fixed $r$ and $s$ in the left-hand side of  \eqref{ss} for  $t = 0$ leads to  \begin{align*} s\left( (A\sharp_{r,0}B)(A\sharp_{s, 1 - 0}B) \right)^2 &= s\left( A^{\frac{r}{2}} (A^{-\frac{1}{2}} B A^{-\frac{1}{2}})^0 A^{\frac{r + s}{2}} (A^{-\frac{1}{2}} B A^{-\frac{1}{2}})^{1 - 0} A^{\frac{s}{2}} \right)^2\\
&= s\left( A^{\frac{2r + s}{2}} A^{-\frac{1}{2}} B A^{-\frac{1}{2}} A^{\frac{s}{2}} \right)^2\\
&= s\left( A^{\frac{2r + s - 1}{2}} B A^{\frac{s - 1}{2}} \right)^2\\
&= \lambda\left( A^{2r + s - 1} B A^{s - 1} B\right)\\
&\succ_{log} \lambda( A^{2r + 2s - 2} B^2 ) \ \ \ \ \ \text{(using Lemma 2.1)}\\
&= s\left( A^{r + s - 1} B \right)^2.\end{align*}
\end{example}

In order to prove our main result, the next lemma is needed and it is known as the Furuta inequality with negative powers, whose proof can be found in \cite{Tan}.

\begin{lem} Let $X, Y$ be two invertible matrices satisfying $0 < Y \leq X$. Let $0 < p' \leq 1, 0 < q' \leq 1$ and $-1 \leq r' < 0.$ Then, it holds that $$X^{\frac{p' + r'}{q'}} \geq \left(X^{\frac{r'}{2}} Y^{p'} X^{\frac{r'}{2}}\right)^{\frac{1}{q'}}$$ as long as the real numbers $p', r'$ and $q'$ satisfy \begin{equation} -r'(1 - q') \leq p'\leq q' - r'(1 - q'),\label{C1}\end{equation} and one of the following two conditions:
\begin{equation} \frac{1}{2}\leq q' \leq 1\label{C2}\end{equation} or \begin{equation} 0\leq q'\leq \frac{1}{2} \ \ \text{and} \ \ \frac{-r'(1 - q') - q'}{1 - 2q'} \leq p' \leq \frac{-r'(1 - q')}{1 - 2q'}.\label{C3}\end{equation}
\end{lem}

Now we are in a position to present our first result which generalizes inequality \eqref{2} and gives a further generalization of inequality \eqref{ss} precisely
in all the cases where it has been proven valid.

\begin{thm} Let $A$ and $B$ be two positive semi-definite matrices. Then, for all $1\leq p \leq 2$, $r,s \in \mathbb{R}$ with similar signs and $t \in \mathbb{R}$ such that $\frac{rp - r}{(r + s)p}\leq t\leq \frac{rp + s}{(r + s)p}$, we have $$\lambda((A\sharp_{r,t} B)^p (A\sharp_{s,1-t}B)^p) \prec_{log} \lambda(A^{r + s - 1}B)^p.$$
\end{thm}

\begin{proof} Without loss of generality, we shall assume that $A$ and $B$ are positive definite matrices as the general case can be then obtained by a continuity argument. Let $r, s \in \mathbb{R}$ with similar signs and $t \in \mathbb{R}$ such that $\frac{rp - r}{(r + s)p}\leq t\leq \frac{rp + s}{(r + s)p}$. Note that in order to finish the proof in this case, it is enough to show that for all $A, B > 0$ the following is true \begin{equation} A^{\frac{r + s - 1}{2}} B A^{\frac{r + s - 1}{2}} \leq I_n \Rightarrow (A\sharp_{s,1 - t} B)^{\frac{p}{2}} (A\sharp_{r,t} B)^p (A\sharp_{s,1 - t}B)^{\frac{p}{2}} \leq I_n.\label{implication}\end{equation}

Now assume that $A^{\frac{r + s - 1}{2}} B A^{\frac{r + s - 1}{2}} \leq I_n$. Then clearly, this is equivalent to \begin{equation} 0< A^{-\frac{1}{2}} B A^{-\frac{1}{2}} \leq
 A^{-(r + s)}.\label{qq}\end{equation}

Applying Lemma 2.2 on \eqref{qq} with $X = A^{-(r + s)}$ and $Y = A^{-\frac{1}{2}} B A^{-\frac{1}{2}}$, we then obtain  \begin{equation}\left(A^{-(r + s)}\right)^{\frac{p' + r'}{q'}} \geq \left[\left(A^{-(r + s)}\right)^{\frac{r'}{2}} \left(A^{-\frac{1}{2}} B A^{-\frac{1}{2}}\right)^{p'} \left(A^{-(r + s)}\right)^{\frac{r'}{2}}\right]^{\frac{1}{q'}}.\label{eq}\end{equation}

Now taking $p' = t$, $r' = -\frac{r}{r + s}$ and $q' = \frac{1}{p}$ in \eqref{eq} for which conditions \eqref{C1} and \eqref{C2} are satisfied, we get $$A^{rp - (r + s)pt} \geq \left(A^{\frac{r}{2}} (A^{-\frac{1}{2}} B A^{-\frac{1}{2}})^t A^{\frac{r}{2}}\right)^p$$ which is the same as  \begin{equation} (A\sharp_{r,t}B)^{p} \leq A^{rp - (r + s)pt}.  \label{3} \end{equation}

Again, replacing this time  $p'$ with $1 - t$, $r'$ with $-\frac{s}{r + s}$ and $q'$ with $\frac{1}{p}$ in \eqref{eq} where  conditions \eqref{C1} and \eqref{C2} are also satisfied,  yields  $$A^{(r + s)pt - rp} \geq \left(A^{\frac{s}{2}} (A^{-\frac{1}{2}} B A^{-\frac{1}{2}})^{1 - t} A^{\frac{s}{2}}\right)^p = (A\sharp_{s,1 - t} B)^p,$$ which is equivalent to \begin{equation} A^{rp - (r + s)pt} \leq (A\sharp_{s,1-t}B)^{-p}.  \label{4}\end{equation}

Therefore, we can write \begin{align*} (A\sharp_{s,1-t}B)^{\frac{p}{2}} (A\sharp_{r,t}B)^p (A\sharp_{s,1-t}B)^{\frac{p}{2}} &\leq (A\sharp_{s,1-t}B)^{\frac{p}{2}} A^{rp - (r + s)pt} (A\sharp_{s,1-t}B)^{\frac{p}{2}} \ \ \ \ \ \ \ \ \ \  \text{(using \eqref{3})}\\
&\leq (A\sharp_{s,1-t}B)^{\frac{p}{2}}(A\sharp_{s,1-t}B)^{-p}(A\sharp_{s,1-t}B)^{\frac{p}{2}} \ \ \ \ \ \ \ \text{(using \eqref{4})}\\
&= I_n.\end{align*}

Let $a = \lambda_1(A^{r + s - 1}B)^p > 0$. We observe that $\left(A^{\frac{r + s - 1}{2}} B A^{\frac{r + s - 1}{2}}\right)^p \leq a I_n$, which is equivalent to $$\left(\frac{1}{a^{\frac{1}{2p(r + s -1)}}}A\right)^{\frac{r + s - 1}{2}} \left(\frac{1}{a^{\frac{1}{2p}}} B\right) \left(\frac{1}{a^{\frac{1}{2p(r + s -1)}}}A\right)^{\frac{r + s - 1}{2}} \leq I_n.$$

For simplicity, we will replace $\left(\frac{1}{a^{\frac{1}{2p(r + s -1)}}} A\right)$ and $\left(\frac{1}{a^{\frac{1}{2p}}} B\right)$ with $A'$ and $B'$, respectively, then $$(A')^{\frac{r + s - 1}{2}} (B') (A')^{\frac{r + s - 1}{2}} \leq I_n$$  which using \eqref{implication} implies that $$(A'\sharp_{s,1-t} B')^{\frac{p}{2}} (A'\sharp_{r,t} B')^p (A'\sharp_{s,1-t} B')^{\frac{p}{2}} \leq I_n.$$

So, $$\frac{1}{a}\left[(A\sharp_{s,1-t} B)^{\frac{p}{2}} (A\sharp_{r,t} B)^p (A\sharp_{s,1-t} B)^{\frac{p}{2}}\right] \leq I_n.$$

Hence, $$(A\sharp_{s,1-t} B)^{\frac{p}{2}} (A\sharp_{r,t} B)^p (A\sharp_{s,1-t} B)^{\frac{p}{2}} \leq \lambda_1(A^{r + s - 1} B)^p I_n.$$

Thus, for all $1\leq p \leq 2$ and for all  $r,s \in \mathbb{R}$ with similar signs and such that $\frac{rp - r}{(r + s)p}\leq t\leq \frac{rp + s}{(r + s)p}$, we have that $$\lambda_1\left((A\sharp_{r,t} B)^p (A\sharp_{s,1-t}B)^p\right) \leq \lambda_1\left(A^{r + s - 1}B\right)^p.$$

Finally, as usual by an anti-symmetric tensor product argument the proof is achieved.
\end{proof}

It is worthy to note here that as a result of Theorem 2.1 and Araki-Lieb-Thirring inequality (see, e.g. \cite{con}), we can conclude that for all $1\leq p\leq 2$ and $\frac{rp - r}{(r + s)p} \leq t\leq \frac{rp + s}{(r + s)p}$, the following inequality holds

\begin{equation}\lambda((A\sharp_{r,t} B)^p (A\sharp_{s,1-t}B)^p) \prec_{log} \lambda(A^{r + s - 1}B)^p \prec_{log} \lambda(A^{p(r + s - 1)}B^p).\label{zz}\end{equation}

Clearly, the particular cases  $p = 1$ and $p = 2$ in \eqref{zz} correspond to R. Lemos and G. Soares results.

The following corollary is a result of Theorem 2.1 for the case $r =  s = 1$.

\begin{cor} Let $A$ and $B$ be two positive semi-definite matrices. Then for all $1\leq p \leq 2$ $$\lambda((A\sharp_{t} B)^p (A\sharp_{1-t}B)^p) \prec_{log} \lambda(AB)^p, \hspace{1cm} \frac{p - 1}{2p}\leq t\leq \frac{p + 1}{2p}.$$
\end{cor}

The proof of the next theorem can be done in a similar fashion as that of Theorem 2.1 with only minor changes; namely by making
 use  of condition \eqref{C3} instead of \eqref{C2} in Lemma 2.2.

\begin{thm} Let $A$ and $B$ be two positive semi-definite matrices. Then for all $p \geq 2$, it holds that $$\lambda((A\sharp_{t} B)^p (A\sharp_{1-t}B)^p) \prec_{log} \lambda(AB)^p, \hspace{1cm} \frac{p - 1}{2p}\leq t\leq \frac{p + 1}{2p}.$$
\end{thm}

Taking into account Example 2.1 and in view of Theorem 2.1 as well as Theorem 2.2, we conclude this section with  the following more general conjecture.

\begin{con} If $A, B \geq 0$, $0\leq t\leq 1$, $p \geq 1$ and ($r, s\geq 1$ or $r, s \leq 0$),  then $$\lambda((A\sharp_{r,t}B)^p(A\sharp_{s, 1 - t}B)^p) \prec_{log} \lambda( A^{p(r + s - 1)}B^p).$$
\end{con}

\section{Results related to Conjecture 1.1}
This section deals with some results related to Conjecture 1.1. In particular, we shall present a new proof of a refinement of inequality \eqref{1} proved by R. Lemos and G. Soares \cite[Corollary 7.2]{Lemos2} (it worth mentioning here  that Corollary 7.2 in \cite{Lemos2} deals with a more general result).

We will start with the following lemma which is a great tool for showing majorization inequalities and whose proof can be found in \cite[p. 352]{Zhang}.

\begin{lem} Let $M = \left[ \begin{array}{cc}
                     A & B\\
                     B^* & C\\
                    \end{array} \right] \geq 0.$ Then $$\vert \lambda(B)\vert \prec_{log} s(B) \prec_{w log} \lambda(A)^{\frac{1}{2}} \circ \lambda(C)^{\frac{1}{2}},$$  where $\circ $ represents the Hadamard product.\end{lem}

For the sake of simplicity, throughout this section we shall use  the following notation $$X = A^{-\frac{1}{2}} B A^{-\frac{1}{2}},$$  whenever $A$ and $B$ are
are both fixed positive semi-definite matrices in $M_n$.

\begin{thm} Let $A$ and $B$ be two positive semi-definite matrices. Then $$s\left(A^{\frac{1}{2}}(A\sharp B)B^{ \frac{1}{2}}\right) \prec_{log} \lambda(AB) \prec_{log} s(AB).$$ \end{thm}

\begin{proof} If $M  = \left[ \begin{array}{cc}
                     B^{\frac{1}{2}} A B^{\frac{1}{2}} & B^{\frac{1}{2}} (A\sharp B) A^{\frac{1}{2}}\\
                     A^{\frac{1}{2}}  (A\sharp B) B^{\frac{1}{2}} & A^{\frac{1}{2}} B A^{\frac{1}{2}}\\
                    \end{array} \right],$ then clearly $M$ can be rewritten as $$M =  \left[ \begin{array}{cc}
                     A^{\frac{1}{2}}B^{\frac{1}{2}} & X^{\frac{1}{2}} A
                    \end{array} \right]^*\left[ \begin{array}{cc}
                     A^{\frac{1}{2}}B^{\frac{1}{2}} & X^{\frac{1}{2}}A \end{array} \right] \geq 0.$$

Taking into account that the determinants of $A^{\frac{1}{2}} (A\sharp B) B^{\frac{1}{2}}$ and $AB$ are equal, then applying Lemma 3.1 on the $2\times 2$ block matrix $M$  gives
$$s\left(A^{\frac{1}{2}} (A\sharp B) B^{\frac{1}{2}}\right) \prec_{ log} \lambda\left(B^{\frac{1}{2}} A B^{\frac{1}{2}} \right)^{\frac{1}{2}} \circ \lambda\left(A^{\frac{1}{2}} B A^{\frac{1}{2}}\right)^{\frac{1}{2}} = \lambda(AB) \prec_{log} s(AB).$$
\end{proof}

\begin{remark} Motivated by Theorem 3.1, one naturally would ask whether the following generalized inequality is true $$s(A^t(A\sharp_tB) B^{1-t}) \prec_{log}\lambda(AB) \hspace{1cm} 0\leq t\leq 1.$$ Unfortunately it is not, since for $t = 0$ we have $$s(A^0(A\sharp_0B) B^{1 - 0}) = s(AB) \succ_{log} \lambda(AB).$$   \end{remark}

The next theorem shows that replacing the singular values with eigenvalues in \eqref{con} is true, however Conjecture 1.1 is still an open problem.

\begin{thm} Let $A$ and $B$ be two positive semi-definite matrices. Then $$\left\vert \lambda\left(A^t (A\sharp_{t}B) B^{1-t}\right)\right\vert \prec_{log} \lambda(AB), \hspace{1cm} 0\leq t\leq 1.$$
\end{thm}

\begin{proof} As usual, we shall assume that $A$ is a positive definite matrix, and the general case can then be deduced by a continuity argument. In addition, we shall divide the proof into two cases.\begin{enumerate}
    \item[\underline{Case 1:}]  If $0\leq t\leq \frac{1}{2}$, then let $$M = \left[ \begin{array}{cc}
                     B^{\frac{1}{2}} A B^{\frac{1}{2}} & B^{\frac{1}{2}}  A^t (A\sharp_{t}B) B^{\frac{1}{2} - t}\\
                     B^{\frac{1}{2} - t} (A\sharp_{t}B) A^t B^{\frac{1}{2}} & B^{\frac{1}{2} - t} (A\sharp_{t}B) A^{-\frac{1}{2}} A^{2t} A^{-\frac{1}{2}} (A\sharp_{t}B) B^{\frac{1}{2} - t}\\
                    \end{array} \right],$$ which can be rewritten in terms of $X = A^{-\frac{1}{2}} B A^{-\frac{1}{2}}$ as $$M = \left[ \begin{array}{cc}
                     B^{\frac{1}{2}} A B^{\frac{1}{2}} & B^{\frac   {1}{2}} A^{\frac{1}{2}} A^t X^t  A^{\frac{1}{2}} B^{\frac{1}{2} - t}\\
                     B^{\frac{1}{2} - t} A^{\frac{1}{2}} X^t  A^t A^{\frac{1}{2}} B^{\frac{1}{2}} & B^{\frac{1}{2} - t} A^{\frac{1}{2}} X^t A^{2t} X^t A^{\frac{1}{2}} B^{\frac{1}{2} - t}\\
                    \end{array} \right].$$

Moreover, $M$ can be also rewritten as $Z^*Z \geq 0$, where $$Z =  \left[ \begin{array}{cc}
A^{\frac{1}{2}}B^{\frac{1}{2}} &  A^t X^t A^{\frac{1}{2}} B^{\frac{1}{2} - t}
\end{array} \right], $$ so that $M\geq 0$.

Now, applying Lemma 3.1 on the matrix $M$ gives
\begin{equation} \left\vert\lambda\left(A^t (A\sharp_{t}B) B^{1 - t}\right)\right\vert \prec_{w log} \lambda\left(B^{\frac{1}{2}} A B^{\frac{1}{2}}\right)^{\frac{1}{2}}\circ \lambda\left(B^{\frac{1}{2} - t} A^{\frac{1}{2}} X^t A^{2t} X^t A^{\frac{1}{2}} B^{\frac{1}{2}
- t}\right)^{\frac{1}{2}} .\label{7}\end{equation}

Next, our goal is to prove that for all $0\leq t\leq \frac{1}{2}$, $$\lambda\left(B^{\frac{1}{2} - t} A^{\frac{1}{2}} X^t A^{2t} X^t A^{\frac{1}{2}} B^{\frac{1}{2} - t}\right) \prec_{log} \lambda\left(AB\right),$$ for which, by using an anti-symmetric tensor product argument,  is sufficient to show that for all $0\leq t\leq \frac{1}{2}$, $$\lambda_1\left(B^{\frac{1}{2} - t} A^{\frac{1}{2}} X^t A^{2t} X^t A^{\frac{1}{2}} B^{\frac{1}{2} - t}\right) \leq \lambda_1\left(AB\right).$$

Without loss of generality, we shall assume that $\lambda_1(AB) = 1$. Obviously, proving our claim is now equivalent to showing that $$\lambda_1\left(B^{\frac{1}{2} - t} A^{\frac{1}{2}} X^t A^{2t} X^t A^{\frac{1}{2}} B^{\frac{1}{2} - t}\right) \leq 1.$$

As $\lambda_1(AB) = 1$, then this implies that $$A^{\frac{1}{2}}BA^{\frac{1}{2}} \leq I_n$$ which in turn gives that \begin{equation} B \leq A^{-1}. \label{8}\end{equation}

Applying Lowner-Heinz inequality (see \cite[Theorem 1.1]{Zhan}) on \eqref{8} for $0\leq 1-2t\leq 1$ we obtain  \begin{equation}B^{1 - 2t} \leq A^{2t - 1}.\label{v}\end{equation}

Next, multiplying both sides of \eqref{8} with $A^{-\frac{1}{2}} > 0$ yields $$X \leq A^{-2},$$ and by appealing again to Lowner-Heinz inequality this time  for $0\leq t\leq \frac{1}{2}$, we obtain \begin{equation} X^t \leq A^{-2t}\label{9}\end{equation}

Now, we can write \begin{align*} \lambda_1\left(B^{\frac{1}{2} - t} A^{\frac{1}{2}} X^t A^{2t} X^t A^{\frac{1}{2}} B^{\frac{1}{2} - t}\right) &= \lambda_1\left(A^t X^t A^{\frac{1}{2}} B^{1- 2t} A^{\frac{1}{2}} X^t A^t\right)\\
&\leq \lambda_1\left(A^t X^t A^{2t} X^t A^t\right)\ \ \ \ \ \text{(using \eqref{v})}\\
&= \lambda_1\left(A^{2t} X^t A^{2t} X^t\right)\\
&= \lambda_1(A^{2t}X^t)^2\\
&\leq \lambda_1(I_n)^2 \hspace{2 cm} \ \ \  \text{(using \eqref{9})}\\
&= 1. \end{align*}

Therefore, \begin{equation} \lambda\left(B^{\frac{1}{2} - t} A^{\frac{1}{2}} X^t A^{2t} X^t A^{\frac{1}{2}} B^{\frac{1}{2} - t}\right) \prec_{log} \lambda\left(AB\right), \hspace{1cm}0\leq t\leq \frac{1}{2}.\label{10}\end{equation}

On the other hand, using \eqref{7} and \eqref{10} gives \begin{align*}\left\vert\lambda\left(A^t (A\sharp_{t}B)  B^{1 - t}\right)\right\vert &\prec_{w log} \lambda\left(B^{\frac{1}{2}} A B^{\frac{1}{2}}\right)^{\frac{1}{2}}\circ \lambda\left(B^{\frac{1}{2} - t} A^{\frac{1}{2}} X^t A^{2t} X^t A^{\frac{1}{2}} B^{\frac{1}{2} - t}\right)^{\frac{1}{2}}\\
&\prec_{w log} \lambda(B^{\frac{1}{2}} A B^{\frac{1}{2}})^{\frac{1}{2}} \circ \lambda(AB)^{\frac{1}{2}}\\
&= \lambda(AB). \end{align*}

Making use of the fact that the determinants of the matrices on the left and on the right are equal, we finally arrive at $\left\vert\lambda\left(A^t (A\sharp_{t}B)  B^{1 - t}\right)\right\vert \prec_{ log } \lambda(AB).$\\

\item[\underline{Case 2:}] Let $\frac{1}{2} \leq t\leq 1$. Replacing $A$ with $B$, $B$ with $A$ and $t$ with $1 - t$ in Case 1, yields $$\left\vert \lambda\left(B^{1 - t} (B\sharp_{1 - t}A) A^{1- (1 - t)}\right)\right\vert \prec_{log} \lambda(BA) \hspace{1cm} 0\leq 1 - t\leq \frac{1}{2}.$$

By noting that $A\sharp_t B = B\sharp_{1-t}A$ and that $\vert \lambda(X)\vert = \vert \lambda(X^*)\vert$ for all $X \in M_n$, we thus have  that for all $\frac{1}{2}\leq t\leq 1$, \begin{align*}\left\vert \lambda\left(A^t (A\sharp_t B) B^{1 - t}\right)\right\vert &= \left\vert \lambda\left(B^{1 - t} (B\sharp_{1 - t} A) A^{t}\right)\right\vert\\
&\prec_{log} \lambda(BA)\\
&= \lambda(AB). \end{align*}
\end{enumerate}
\end{proof}

Our next goal is to show a result that strengthens Conjecture 1.1. For this,  we need the following lemma whose proof can be found in \cite{AGM}.

\begin{lem} Let $A$ be a positive semi-definite matrix and $B$ be a Hermitian matrix. Then for all $p, q \geq 0$, $$\lambda(A^p B A^{q} B) \prec_{log} \lambda(A^{p + q} B^2).$$ \end{lem}

\begin{thm} Let $A$ and $B$ be two positive definite matrices. Then,
\begin{enumerate}
\item $s(A^t (A\sharp_tB) B^{1 - t}) \prec_{log} s(A^{\frac{3}{2}}BA^{-\frac{1}{2}})$ for $\frac{1}{2} \leq t \leq 1$,
\item $s(A^t (A\sharp_tB) B^{1 - t}) \prec_{log} s(B^{\frac{3}{2}}AB^{-\frac{1}{2}})$ for $0 \leq t \leq \frac{1}{2}.$
\end{enumerate}
\end{thm}

\begin{proof} \begin{enumerate}
\item  Let $\frac{1}{2}\leq t\leq 1$. Without loss of generality, we shall again assume that $\lambda_1(A^{\frac{3}{2}}BA^{-1}BA^{\frac{3}{2}}) = 1$. As earlier, our claim now amounts to proving that  $$\lambda_1(A^{t} (A\sharp_tB) B^{2(1 - t)} (A\sharp_tB) A^{t}) \leq 1.$$

Noting that $\lambda_1(A^{\frac{3}{2}}BA^{-1}BA^{\frac{3}{2}}) = \lambda_1(A^2 X^2 A^2)$,  which clearly implies that $X^{2} \leq A^{-4}$, we then obtain \begin{equation} X^{2t} \leq A^{-4t} \label{k}\end{equation} and \begin{equation} B^{2(1 - t)} \leq A^{-2(1 - t)}.\label{l}\end{equation}

Now, we can write \begin{align*} &\lambda_1( A^{t} (A\sharp_tB) B^{2(1-t)} (A\sharp_tB) A^{t} )\\
                    &\leq \lambda_1( A^{t + \frac{1}{2}} X^t A^{\frac{1}{2}} A^{-2(1 - t)} A^{\frac{1}{2}} X^t A^{t + \frac{1}{2}} ) \ \ \ \ \ \text{(using \eqref{l})}\\
                    &= \lambda_1( A^{t + \frac{1}{2}} X^t A^{2t - 1} X^t A^{t + \frac{1}{2}} )\\
                    &= \lambda_1( A^{2t + 1} X^t A^{2t - 1} X^t )\\
                    &\leq \lambda_1( A^{4t} X^{2t} ) \ \  \hspace{4 cm} \text{(using Lemma 3.2)}\\
                    &= \lambda_1(I_n) \ \ \ \ \ \ \ \ \ \ \   \hspace{4 cm}  \text{(using \eqref{k})}.\\
                    &= 1.\end{align*}

Hence, \begin{equation}\lambda_1( A^{t} (A\sharp_tB) B^{2(1 - t)} (A\sharp_tB) A^{t} ) \leq \lambda_1( A^3BA^{-1}B ), \hspace{1cm}  \frac{1}{2} \leq t \leq 1.\label{x}\end{equation}

Applying square root on both sides of Inequality \eqref{x} gives $$s_1( A^{t} (A\sharp_tB) B^{1 - t} ) \leq s_1(A^{\frac{3}{2}}BA^{-\frac{1}{2}}), \hspace{2.9cm}  \frac{1}{2} \leq t \leq 1.$$

Therefore, with an anti-symmetric tensor product argument the proof can be achieved.

\item If $0 \leq t \leq \frac{1}{2},$ then as was done earlier, replacing $A$ with $B$, $B$ with $A$ and $t$ with $1 - t$ in the preceding case gives $$s(B^{1 - t} (B\sharp_{1 - t}A) A^{t}) \prec_{log} s(B^{\frac{3}{2}}AB^{-\frac{1}{2}}).$$

Finally, to complete the proof it is enough to notice that $B\sharp_{1 - t} A = A\sharp_t B$.
\end{enumerate}
\end{proof}

Before we close this section, it is worthy to note  that in view of  Lemma 2.1, we conclude that for all $A, B > 0$, $$s(AB) \prec_{log} s(A^{\frac{3}{2}}BA^{-\frac{1}{2}}) \ \ \text{and} \ \ s(AB) \prec_{log} s(B^{\frac{3}{2}}AB^{-\frac{1}{2}}).$$

\section{A complement of a norm inequality}

In \cite{Bhatia2}, R. Bhatia, Y. Lim and T. Yamazaki showed the following norm inequality for all $A, B \geq 0$,  \begin{equation}\vert\vert A + B + 2(A\sharp B)\vert\vert_p \leq \vert\vert A + B + A^{\frac{1}{2}}B^{\frac{1}{2}} + B^{\frac{1}{2}}A^{\frac{1}{2}}\vert\vert_p \hspace{1cm} p = 1, 2, \infty.\label{11}\end{equation}

Then they conjectured that it is true for every unitarily invariant norm, which was recently proved by T.H. Dinh, R. Dumitru and J.A. Franco \cite{Dinh} for the Schatten p-norms with $1\leq p\leq \infty$.\\

First, we shall borrow from \cite{Lin} the following definition for  positive definite matrices $A$ and $B$, $$A \natural\natural  B = A^{\frac{1}{2}}(B^{\frac{1}{2}}A^{-1}B^{\frac{1}{2}})^{\frac{1}{2}}A^{\frac{1}{2}}.$$

In addition, we shall need the following basic result.
\begin{prop} Let $M:=\left[\begin{array}{cc}
                     A & B\\
                     B & A\\
                    \end{array} \right]$ where $A$ and $B$ are in $M_n$. Then, the eigenvalues of $M$ are the union of the eigenvalues of $A+B$ and $A-B$.

\end{prop}

\begin{thm} Let $A$, $B$ and $C$ be positive definite matrices of order $n$ such that $C \geq A + B$. Then for all $j = 1, 2, \dots, n$,  $$s_j(C + A^{\frac{1}{2}}B^{\frac{1}{2}} + B^{\frac{1}{2}}A^{\frac{1}{2}})\leq s_j(C + A\sharp B + A\natural\natural B).$$ \end{thm}

\begin{proof} Let $M = \left[ \begin{array}{cc}
                     A\sharp B & B^{\frac{1}{2}}A^{\frac{1}{2}}\\
                     A^{\frac{1}{2}}B^{\frac{1}{2}} & A\natural\natural B\\
                    \end{array} \right] \geq 0.$ Then $N = \left[ \begin{array}{cc}
                     A\natural\natural B & A^{\frac{1}{2}}B^{\frac{1}{2}}\\
                     B^{\frac{1}{2}}A^{\frac{1}{2}} & A\sharp B\\
                    \end{array} \right] \geq 0.$

Hence, $$T = M + N +  \left[\begin{array}{cc}
                     C & C\\
                     C & C\\
                    \end{array} \right] = \left[ \begin{array}{cc}
                     C + A\sharp B + A\natural\natural B & C + A^{\frac{1}{2}}B^{\frac{1}{2}} + B^{\frac{1}{2}}A^{\frac{1}{2}}\\
                     C + A^{\frac{1}{2}}B^{\frac{1}{2}} + B^{\frac{1}{2}}A^{\frac{1}{2}} & C + A\sharp B + A\natural\natural B\\
                    \end{array} \right] \geq 0.$$

From the preceding proposition, we conclude that  $$C + A^{\frac{1}{2}}B^{\frac{1}{2}} + B^{\frac{1}{2}}A^{\frac{1}{2}} \leq C + A\sharp B + A\natural\natural B \ \ \text{and} \ \ -(C + A^{\frac{1}{2}}B^{\frac{1}{2}} + B^{\frac{1}{2}}A^{\frac{1}{2}}) \leq C + A\sharp B + A\natural\natural B.$$

Thus, $$\lambda_j(C + A^{\frac{1}{2}}B^{\frac{1}{2}} + B^{\frac{1}{2}}A^{\frac{1}{2}}) \leq \lambda_j(C + A\sharp B + A\natural\natural B).$$

Consequently, from the fact that $C + A^{\frac{1}{2}}B^{\frac{1}{2}} + B^{\frac{1}{2}}A^{\frac{1}{2}} \geq A + B + A^{\frac{1}{2}}B^{\frac{1}{2}} + B^{\frac{1}{2}}A^{\frac{1}{2}} \geq 0$ we get the result.

\end{proof}

The following corollary can be considered as a complement of inequality \eqref{11} when replacing $C$ by $A + B$ in Theorem 4.1.

\begin{cor} Let $A$ and $B$ be two positive definite matrices of order $n$. Then for all unitarily invariant norms,  $$\vert\vert\vert A + B + A^{\frac{1}{2}}B^{\frac{1}{2}} + B^{\frac{1}{2}}A^{\frac{1}{2}}\vert\vert\vert \leq \vert\vert\vert A + B + A\sharp B + A\natural\natural B\vert\vert\vert.$$ \end{cor}

Next, we present a further complement of the preceding corollary.

\begin{thm} Let $A$ and $B$ be two positive definite matrices of order $n$. Then, for $p = 1, 2$,  $$\vert\vert A + B + A\sharp B + A\natural\natural B\vert\vert_p \leq \vert\vert A + B + 2(A\natural\natural B)\vert\vert_p.$$ \end{thm}

Before presenting the proof, we shall consider the following log-majorization relations. The first one is due to M. Lin \cite{Lin}.

\begin{lem} Let $A$ and $B$ be two positive definite matrices. Then, $$\lambda(A\sharp B) \prec_{log} \lambda(A\natural\natural B).$$
\end{lem}

\begin{lem} Let $A$ and $B$ be two positive definite matrices. Then, $$\lambda(A (A\sharp B)) \prec_{log} \lambda(A (A\natural\natural B)).$$
\end{lem}

\begin{proof} Recalling first an inequality from \cite{AG} which can be stated for $A, B \geq 0$ as follows \begin{equation}\lambda(A^{\frac{k}{2}}(A^{-\frac{1}{2}} B A^{-\frac{1}{2}})^t A^{\frac{k}{2}}) \prec_{log} \lambda(A^{k - t}B^t) \hspace{1cm} 0\leq t\leq 1 \ \ \text{and} \ \ k\geq t,\label{w}\end{equation} and another one from \cite{GAMA} which is a complement of \eqref{w} and can be stated as  $$\lambda(A^{k - t}B^t) \prec_{log} \lambda(A^{\frac{k}{2}}(B^{\frac{1}{2}} A^{-1} B^{\frac{1}{2}})^t A^{\frac{k}{2}}) \hspace{1cm} 0\leq t\leq 1 \ \ \text{and} \ \  k\leq 2.$$

Now, taking $k = 2$ and $t = \frac{1}{2}$ in both inequalities, we obtain $$\lambda(A (A^{-\frac{1}{2}} B A^{-\frac{1}{2}})^{\frac{1}{2}} A) \prec_{log} \lambda(A^{\frac{3}{2}}B^{\frac{1}{2}}) \prec_{log} \lambda(A (B^{\frac{1}{2}} A^{-1} B^{\frac{1}{2}})^{\frac{1}{2}} A)$$ which gives the desired result.
\end{proof}

\begin{remark} It is worthy to mention here that the log-majorization relation \eqref{w} was first proved by Furuta \cite[Corollary 3.1]{Furuta}, and for $k = 1$ it was obtained earlier  by  J.S. Matharu and J.S. Aujla \cite[Theorem 2.10]{Mat}. Later, for $k = 2$ and $t = \frac{1}{2}$, it was reobtained by Bhatia, Lim and Yamazaki \cite[Theorem 2]{Bhatia2}, and also for $k = 2$ it was reobtained  by  D.T. Hoa \cite[Proposition 2.1]{Hoa}.\end{remark}

\begin{lem} Let $A$ and $B$ be two positive definite matrices. Then $$\lambda(B (A\sharp B)) \prec_{log} \lambda(B (A\natural\natural B)).$$
\end{lem}

\begin{proof} By Schur's complement we know that $$M = \left[ \begin{array}{cc}
                     B^{\frac{1}{2}}A^{\frac{1}{2}}(B^{\frac{1}{2}}A^{-1}B^{\frac{1}{2}})^{\frac{1}{2}} A^{\frac{1}{2}}B^{\frac{1}{2}} & B^{\frac{1}{2}}A^{\frac{1}{2}}B\\
                     BA^{\frac{1}{2}}B^{\frac{1}{2}} & B(B^{-\frac{1}{2}}AB^{-\frac{1}{2}})^{\frac{1}{2}}B\\
                    \end{array} \right] \geq 0.$$

Next, by appealing to Lemma 3.1, we obtain \begin{align*}\lambda_1(B^{\frac{3}{2}}A^{\frac{1}{2}})^2 &\leq \lambda_1(B^{\frac{1}{2}}A^{\frac{1}{2}}(B^{\frac{1}{2}}A^{-1}B^{\frac{1}{2}})^{\frac{1}{2}})A^{\frac{1}{2}}B^{\frac{1}{2}})\cdot\lambda_1(B(B^{-\frac{1}{2}}AB^{-\frac{1}{2}})^{\frac{1}{2}}B)\\
&= \lambda_1(B (A\natural\natural B))\cdot \lambda_1(B (B\sharp A))\\
&\leq \lambda_1(B (A\natural\natural B))\cdot \lambda_1(BB^{\frac{1}{2}}A^{\frac{1}{2}})\ \ \ \ \ \text{(using \eqref{w})}\\
\end{align*} but this  implies that  $$\lambda_1(B (A\sharp B)) \leq \lambda_1(B^{\frac{3}{2}}A^{\frac{1}{2}}) \leq \lambda_1(B (A\natural\natural B)).$$

Finally, we obtain the result with a standard anti-symmetric tensor product argument.
\end{proof}

Now, we are in position to prove Theorem 4.2.\\

\textbf{Proof of Theorem 4.2.} Using Lemma 4.1, it is easy to see that $$Tr(A + B + A\sharp B + A\natural\natural B) \leq Tr(A + B + 2(A\natural\natural B)),$$ so that for $p = 1$, the inequality is satisfied. Now, for $p = 2$, it is enough to show that \begin{align*} &Tr(A\sharp B) \leq Tr(A\natural\natural B);\\
&Tr(A (A\sharp B)) \leq Tr(A (A\natural\natural B));\\
&Tr(B (A\sharp B)) \leq Tr(B (A\natural\natural B)),\end{align*} which, in view of the previous lemmas, these are all true. This completes the proof.  \hspace{1cm} $\Box$ \\

Based on our work in this section, we conclude the paper with the following conjecture.

\begin{con} Let $A$ and $B$ be two positive definite matrices of order $n$. Then for $1\leq p\leq \infty$  $$\vert\vert A + B + A\sharp B + A\natural\natural B\vert\vert_p \leq \vert\vert A + B + 2(A\natural\natural B)\vert\vert_p.$$
\end{con}

Finally, it is worthy to observe that  the validation of the previous conjecture would imply that for all $1\leq p\leq \infty$, $$\vert\vert A + B + 2(A\sharp B)\vert\vert_p \leq \vert\vert A + B + A^{\frac{1}{2}}B^{\frac{1}{2}} + B^{\frac{1}{2}}A^{\frac{1}{2}}\vert\vert_p \leq \vert\vert A + B + 2(A\natural\natural B)\vert\vert_p.$$

\section*{Acknowledgements}
The first 3 authors acknowledge financial support from the Lebanese University research grants program.
\section*{Declaration of Competing Interest}
The authors declare that they have no competing interests.

\end{document}